\documentclass[11pt, a4paper]{article}
\usepackage[cp1251]{inputenc}
\usepackage{amsmath} \usepackage{euscript} \usepackage{marvosym}
\usepackage[dvipsnames]{xcolor}
\usepackage{graphicx}
\usepackage{longtable}
\usepackage{mathrsfs}
\usepackage{amssymb, amscd}
\usepackage{comment}
\begin{document}

\newcounter{bnomer} \newcounter{snomer}
\newcounter{bsnomer}
\setcounter{bnomer}{0}
\renewcommand{\thesnomer}{\thebnomer.\arabic{snomer}}
\renewcommand{\thebsnomer}{\thebnomer.\arabic{bsnomer}}
\renewcommand{\refname}{\begin{center}\large{\textbf{References}}\end{center}}

\setcounter{MaxMatrixCols}{14}

\newcommand\restr[2]{{
  \left.\kern-\nulldelimiterspace 
  #1 
  \right|_{#2} 
}}

\newcommand{\sect}[1]{%
\setcounter{snomer}{0}\setcounter{bsnomer}{0}
\refstepcounter{bnomer}
\par\bigskip\begin{center}\large{\textbf{\arabic{bnomer}. {#1}}}\end{center}}
\newcommand{\sst}[1]{%
\refstepcounter{bsnomer}
\par\bigskip\textbf{\arabic{bnomer}.\arabic{bsnomer}. {#1}}\par}
\newcommand{\defi}[1]{%
\refstepcounter{snomer}
\par\medskip\textbf{Definition \arabic{bnomer}.\arabic{snomer}. }{#1}\par\medskip}
\newcommand{\theo}[2]{%
\refstepcounter{snomer}
\par\textbf{Theorem \arabic{bnomer}.\arabic{snomer}. }{#2} {\emph{#1}}\hspace{\fill}$\square$\par}
\newcommand{\mtheop}[2]{%
\refstepcounter{snomer}
\par\medskip\textbf{Theorem \arabic{bnomer}.\arabic{snomer}. }{\emph{#1}}
\par\textsc{Proof}. {#2}\hspace{\fill}$\square$\par}
\newcommand{\mcorop}[2]{%
\refstepcounter{snomer}
\par\textbf{Corollary \arabic{bnomer}.\arabic{snomer}. }{\emph{#1}}
\par\textsc{Proof}. {#2}\hspace{\fill}$\square$\par}
\newcommand{\mtheo}[1]{%
\refstepcounter{snomer}
\par\medskip\textbf{Theorem \arabic{bnomer}.\arabic{snomer}. }{\emph{#1}}\par\medskip}
\newcommand{\theobn}[1]{%
\par\medskip\textbf{Theorem. }{\emph{#1}}\par\medskip}
\newcommand{\theoc}[2]{%
\refstepcounter{snomer}
\par\medskip\textbf{Theorem \arabic{bnomer}.\arabic{snomer}. }{#1} {\emph{#2}}\par\medskip}
\newcommand{\mlemm}[1]{%
\refstepcounter{snomer}
\par\medskip\textbf{Lemma \arabic{bnomer}.\arabic{snomer}. }{\emph{#1}}\par\medskip}
\newcommand{\mprop}[1]{%
\refstepcounter{snomer}
\par\medskip\textbf{Proposition \arabic{bnomer}.\arabic{snomer}. }{\emph{#1}}\par\medskip}
\newcommand{\theobp}[2]{%
\refstepcounter{snomer}
\par\textbf{Theorem \arabic{bnomer}.\arabic{snomer}. }{#2} {\emph{#1}}\par}
\newcommand{\theop}[2]{%
\refstepcounter{snomer}
\par\textbf{Theorem \arabic{bnomer}.\arabic{snomer}. }{\emph{#1}}
\par\textsc{Proof}. {#2}\hspace{\fill}$\square$\par}
\newcommand{\theosp}[2]{%
\refstepcounter{snomer}
\par\textbf{Theorem \arabic{bnomer}.\arabic{snomer}. }{\emph{#1}}
\par\textsc{Sketch of the proof}. {#2}\hspace{\fill}$\square$\par}
\newcommand{\exam}[1]{%
\refstepcounter{snomer}
\par\medskip\textbf{Example \arabic{bnomer}.\arabic{snomer}. }{#1}\par\medskip}
\newcommand{\deno}[1]{%
\refstepcounter{snomer}
\par\textbf{Notation \arabic{bnomer}.\arabic{snomer}. }{#1}\par}
\newcommand{\lemm}[1]{%
\refstepcounter{snomer}
\par\textbf{Lemma \arabic{bnomer}.\arabic{snomer}. }{\emph{#1}}\hspace{\fill}$\square$\par}
\newcommand{\lemmp}[2]{%
\refstepcounter{snomer}
\par\medskip\textbf{Lemma \arabic{bnomer}.\arabic{snomer}. }{\emph{#1}}
\par\textsc{Proof}. {#2}\hspace{\fill}$\square$\par\medskip}
\newcommand{\coro}[1]{%
\refstepcounter{snomer}
\par\textbf{Corollary \arabic{bnomer}.\arabic{snomer}. }{\emph{#1}}\hspace{\fill}$\square$\par}
\newcommand{\mcoro}[1]{%
\refstepcounter{snomer}
\par\textbf{Corollary \arabic{bnomer}.\arabic{snomer}. }{\emph{#1}}\par\medskip}
\newcommand{\corop}[2]{%
\refstepcounter{snomer}
\par\textbf{Corollary \arabic{bnomer}.\arabic{snomer}. }{\emph{#1}}
\par\textsc{Proof}. {#2}\hspace{\fill}$\square$\par}
\newcommand{\nota}[1]{%
\refstepcounter{snomer}
\par\medskip\textbf{Remark \arabic{bnomer}.\arabic{snomer}. }{#1}\par\medskip}
\newcommand{\cons}[1]{%
\refstepcounter{snomer}
\par\medskip\textbf{Construction \arabic{bnomer}.\arabic{snomer}. }{#1}\par\medskip}
\newcommand{\propp}[2]{%
\refstepcounter{snomer}
\par\medskip\textbf{Proposition \arabic{bnomer}.\arabic{snomer}. }{\emph{#1}}
\par\textsc{Proof}. {#2}\hspace{\fill}$\square$\par\medskip}
\newcommand{\hypo}[1]{%
\refstepcounter{snomer}
\par\medskip\textbf{Conjecture \arabic{bnomer}.\arabic{snomer}. }{\emph{#1}}\par\medskip}
\newcommand{\prop}[1]{%
\refstepcounter{snomer}
\par\textbf{Proposition \arabic{bnomer}.\arabic{snomer}. }{\emph{#1}}\hspace{\fill}$\square$\par}

\newcommand{\proof}[2]{%
\par\medskip\textsc{Proof{#1}}. \hspace{-0.2cm}{#2}\hspace{\fill}$\square$\par\medskip}

\makeatletter
\def\iddots{\mathinner{\mkern1mu\raise\p@
\vbox{\kern7\p@\hbox{.}}\mkern2mu
\raise4\p@\hbox{.}\mkern2mu\raise7\p@\hbox{.}\mkern1mu}}
\makeatother

\newcommand{\okr}[2]{%
\refstepcounter{snomer}
\par\medskip\textbf{{#1} \arabic{bnomer}.\arabic{snomer}. }{\emph{#2}}\par\medskip}

\newcommand{\Ind}[3]{%
\mathrm{Ind}_{#1}^{#2}{#3}}
\newcommand{\Res}[3]{%
\mathrm{Res}_{#1}^{#2}{#3}}
\newcommand{\epsi}{\varepsilon}
\newcommand{\tri}{\triangleleft}
\newcommand{\Supp}[1]{%
\mathrm{Supp}(#1)}
\newcommand{\SSu}[1]{%
\mathrm{SingSupp}(#1)}

\newcommand{\gee}{\geqslant}
\newcommand{\reg}{\mathrm{reg}}
\newcommand{\Dyn}{\mathrm{Dyn}}
\newcommand{\Ann}{\mathrm{Ann}\,}
\newcommand{\Cent}[1]{\mathbin\mathrm{Cent}({#1})}
\newcommand{\PCent}[1]{\mathbin\mathrm{PCent}({#1})}
\newcommand{\Irr}[1]{\mathbin\mathrm{Irr}({#1})}
\newcommand{\Exp}[1]{\mathbin\mathrm{Exp}({#1})}
\newcommand{\empr}[2]{[-{#1},{#1}]\times[-{#2},{#2}]}
\newcommand{\sreg}{\mathrm{sreg}}
\newcommand{\ilm}{\varinjlim}
\newcommand{\wdth}{\mathrm{wd}}
\newcommand{\plm}{\varprojlim}
\newcommand{\codim}{\mathrm{codim}\,}
\newcommand{\GKdim}{\mathrm{GKdim}\,}
\newcommand{\SAut}{\mathrm{SAut}}
\newcommand{\Aut}{\mathrm{Aut}}
\newcommand{\Spec}{\mathrm{Spec}\,}
\newcommand{\chara}{\mathrm{char}\,}
\newcommand{\rk}{\mathrm{rk}\,}
\newcommand{\chr}{\mathrm{ch}\,}
\newcommand{\Ker}{\mathrm{Ker}\,}
\newcommand{\id}{\mathrm{id}}
\newcommand{\Ad}{\mathrm{Ad}}
\newcommand{\Gh}{\mathrm{Gh}}
\newcommand{\col}{\mathrm{col}}
\newcommand{\row}{\mathrm{row}}
\newcommand{\high}{\mathrm{high}}
\newcommand{\low}{\mathrm{low}}
\newcommand{\pho}{\hphantom{\quad}\vphantom{\mid}}
\newcommand{\fho}[1]{\vphantom{\mid}\setbox0\hbox{00}\hbox to \wd0{\hss\ensuremath{#1}\hss}}
\newcommand{\wt}{\widetilde}
\newcommand{\wh}{\widehat}
\newcommand{\ad}[1]{\mathrm{ad}_{#1}}
\newcommand{\tr}{\mathrm{tr}\,}
\newcommand{\GL}{\mathrm{GL}}
\newcommand{\SL}{\mathrm{SL}}
\newcommand{\SO}{\mathrm{SO}}
\newcommand{\Or}{\mathrm{O}}
\newcommand{\Sp}{\mathrm{Sp}}
\newcommand{\Sa}{\mathrm{S}}
\newcommand{\Ua}{\mathrm{U}}
\newcommand{\Andre}{\mathrm{Andre}}
\newcommand{\Aord}{\mathrm{Aord}}
\newcommand{\Mat}{\mathrm{Mat}}
\newcommand{\Stab}{\mathrm{Stab}}
\newcommand{\htt}{\mathfrak{h}}
\newcommand{\spt}{\mathfrak{sp}}
\newcommand{\slt}{\mathfrak{sl}}
\newcommand{\sot}{\mathfrak{so}}
\newcommand{\Ga}{\mathbb{G}_a}

\newcommand{\vfi}{\varphi}
\newcommand{\aad}{\mathrm{ad}}
\newcommand{\vpi}{\varpi}
\newcommand{\teta}{\vartheta}
\newcommand{\Bfi}{\Phi}
\newcommand{\Fp}{\mathbb{F}}
\newcommand{\Rp}{\mathbb{R}}
\newcommand{\Zp}{\mathbb{Z}}
\newcommand{\Cp}{\mathbb{C}}
\newcommand{\Ap}{\mathbb{A}}
\newcommand{\Pp}{\mathbb{P}}
\newcommand{\Kp}{\mathbb{K}}
\newcommand{\Np}{\mathbb{N}}
\newcommand{\Vp}{\mathbb{V}}
\newcommand{\ut}{\mathfrak{u}}
\newcommand{\at}{\mathfrak{a}}
\newcommand{\glt}{\mathfrak{gl}}
\newcommand{\hei}{\mathfrak{hei}}
\newcommand{\nt}{\mathfrak{n}}
\newcommand{\kt}{\mathfrak{k}}
\newcommand{\mt}{\mathfrak{m}}
\newcommand{\rt}{\mathfrak{r}}
\newcommand{\rad}{\mathfrak{rad}}
\newcommand{\bt}{\mathfrak{b}}
\newcommand{\unt}{\underline{\mathfrak{n}}}
\newcommand{\gt}{\mathfrak{g}}
\newcommand{\vt}{\mathfrak{v}}
\newcommand{\pt}{\mathfrak{p}}
\newcommand{\Xt}{\mathfrak{X}}
\newcommand{\Po}{\mathcal{P}}
\newcommand{\PV}{\mathcal{PV}}
\newcommand{\Uo}{\EuScript{U}}
\newcommand{\Fo}{\EuScript{F}}
\newcommand{\Do}{\EuScript{D}}
\newcommand{\Eo}{\EuScript{E}}
\newcommand{\Jo}{\EuScript{J}}
\newcommand{\Iu}{\mathcal{I}}
\newcommand{\Mo}{\mathcal{M}}
\newcommand{\Nu}{\mathcal{N}}
\newcommand{\Ro}{\mathcal{R}}
\newcommand{\Co}{\mathcal{C}}
\newcommand{\Ko}{\mathcal{K}}
\newcommand{\So}{\mathcal{S}}
\newcommand{\Lo}{\mathcal{L}}
\newcommand{\Ou}{\mathcal{O}}
\newcommand{\Uu}{\mathcal{U}}
\newcommand{\Tu}{\mathcal{T}}
\newcommand{\Au}{\mathcal{A}}
\newcommand{\Vu}{\mathcal{V}}
\newcommand{\Du}{\mathcal{D}}
\newcommand{\Bu}{\mathcal{B}}
\newcommand{\Sy}{\mathcal{Z}}
\newcommand{\Sb}{\mathcal{F}}
\newcommand{\Gr}{\mathcal{G}}
\newcommand{\Xu}{\mathcal{X}}
\newcommand{\Op}{\mathbb{O}}
\newcommand{\chv}{\mathrm{chv}}
\newcommand{\rtc}[1]{C_{#1}^{\mathrm{red}}}
\newcommand{\dd}{\partial}

\author{Mikhail Ignatev\and Timofey Vilkin}
\date{}
\title{On flexibility of trinomial varieties}\maketitle
\begin{abstract} Trinomial varieties are affine varieties given by a system of equations consisting of polynomials with three terms. Such varieties are total coordinate spaces
of normal varieties with torus action of complexity one. For an affine variety $X$ we consider the subgroup $\SAut(X)$ of the automorphism group generated by all algebraic
subgroups isomorphic to the additive group of the ground field. By definition, an affine variety is flexible if $\SAut(X)$ acts transitively on its regular locus. Gaifullin proved a sufficient condition for a trinomial hypersurface to be flexible. We give a generalization of his results, proving a sufficient condition to be flexible for an arbitrary trinomial variety.

\end{abstract}

\sect{Introduction}

\let\thefootnote\relax\footnote{The article was prepared within the framework of the project ``International Academic Cooperation'' HSE University.}

Let $\Kp$ be an algebraically closed field of characteristic zero, $\Ga$ be its additive group, $X$ be an irreducible affine variety over $\Kp$, and $\Aut(X)$ be the group of regular automorphisms of $X$. One can consider its subgroup $\SAut(X)$ of special automorphisms generated by $\Ga$-subgroups. By definition, a $\Ga$-\emph{subgroup} of $\Aut(X)$ is the image of $\Ga$ in $\Aut(X)$ obtained from a regular action of the group $\Ga$ on the variety $X$. Such an action is called a $\Ga$-\emph{action}.

It turned out that $\Ga$-actions are closely related to locally nilpotent derivations of the algebra $\Kp[X]$ of regular functions on $X$. Given a $\Kp$-algebra $A$, a linear operator $\dd\colon A\to A$ is called a \emph{derivation} if it satisfies the Leibniz rule: $$\dd(ab)=a\dd(b)+\dd(a)b\text{ for all }a,b\in A.$$ If for any $a\in A$ there exists a positive integer $n$ such that $\dd^n(a)=0$ then $\dd$ is called {locally nilpotent}~(LND). For any $t\in\Kp$, one can define the exponent $\exp(t\dd) $ of an LND $\dd$. According to \cite[1.5.1]{Freudenburg17}, the mapping $$\dd\mapsto\{\exp(t\dd),~t\in\Kp\}$$ establishes a one-to-one correspondence between LNDs of $\Kp[X]$ and algebraic subgroups of $\Aut(X)$ isomorphic to $\Ga$. 

In this paper we study an important geometric property of the variety $X$ called flexibility. The variety $X$ is called \emph{flexible} if for every regular point $x\in X$ the tangent space $T_xX$ is spanned by tangent vectors to orbits for various $\Ga$-actions. Flexible varieties were investigated in \cite{ArzhantsevFlennerKalimanKutzschebauchZaidenberg13}. Recall that an action of a group $G$ on a set $X$ is called \emph{infinitely transitive} if it is $m$-transitive for each positive integer $m$, while $m$-transitivity means that for each two $m$-tuples of different elements $x_1,\ldots,x_n$ and $y_1,\ldots,y_n$ of $X$ there exists $g\in G$ such that $g\cdot x_i=y_i$ for all $i$. It turned out that flexibility is equivalent to transitivity and, at the same time, to infinite transitivity of the group of special automorphisms on the set of regular points of $X$ if $\dim X \geq 2$, see\break \cite[Theorem~0.1]{ArzhantsevFlennerKalimanKutzschebauchZaidenberg13} for the detail. Flexibility of some interesting classes of affine varieties was studied, e.g., in   \cite{ArzhantsevFlennerKalimanKutzschebauchZaidenberg13, Gaifullin19, GaifullinShafarevich19}. In some sense, flexible varieties have a lot of $\Ga$-actions.

The paper is devoted to the trinomial varieties introduced by Hausen and Wrobel in the paper \cite{HausenWrobel17}. To define them, we need to introduce some notation. Namely, fix an integer $k\geq2$, a non-negative integer $n_0$, and, for each $i\in\{1,\ldots,k\}$, a positive integer $n_i$. We will consider the ring of polynomials in the variables $T_{ij}$, $0\leq i\leq k$, $1\leq j\leq n_i$. For each $i\in\{0,1,2,\ldots,k\}$, fix a tuple $l_i = (l_{i1},\ldots , l_{in_i})$ of positive integers and define the monomial
\begin{equation}T^{l_i}_i=T^{l_{i1}}_{i1}\ldots T^{l_{in_i}}_{in_i}.\label{formula:T_i}
\end{equation}
Here, if $n_0=0$ (and, consequently, $l_0$ is the empty tuple) then we set $T_0^{l_0}=1$.

Also, fix distinct scalars $\lambda_2,\ldots,\lambda_k\in\Kp^{\times}$, where, as usual, $\Kp^{\times}=\Kp\setminus\{0\}$. By definition, a \emph{trinomial variety} is an affine subvariety of the affine space defined by systems of polynomial equations of the form
\begin{equation}
\begin{cases}
\lambda_2T_0^{l_0} + T_1^{l_1}-T_2^{l_2}=0,\\
\lambda_3T_0^{l_0} + T_1^{l_1}-T_3^{l_3}=0,\\
\ldots,\\
\lambda_kT_0^{l_0} + T_1^{l_1}-T_k^{l_k}=0.
\end{cases}\label{formula:trinomial_variety}
\end{equation}
(See also Definition~\ref{def:tr_var} below for an equivalent description of trinomial varieties from \cite[Construction 1.1]{HausenWrobel17}.) In particular, a \emph{trinomial hypersurface} is by definition a trinomial variety defined by a single equation of the form (\ref{formula:trinomial_variety}).

For instance, the group $\SL_2(\Kp)$ of $2\times2$ matrices with determinant 1 is a trinomial hypersurface in the affine space of all $2\times2$ matrices with $n_0=0$, $n_1=n_2=2$, $l_1=l_2=(1,1)$ and $\lambda_2=-1$:
\begin{equation*}
\begin{pmatrix}
T_{11}&T_{21}\\T_{22}&T_{12}
\end{pmatrix}\in\SL_2(\Kp)\text{ if and only if }-1+T_{11}T_{12}-T_{21}T_{22}=0.
\end{equation*}

Recall that an action of an algebraic torus on an algebraic variety has \emph{complexity one} if a generic orbit has codimension one. Note that each trinomial variety admits a regular action of a torus of complexity one. Trinomial varieties are interesting for various reasons. For example, every normal rational variety $X$ with only constant invertible functions, finitely generated divisor class group and an algebraic torus
action of complexity one can be obtained as a quotient of a trinomial variety via action of
a diagonalizable group, see \cite[Corollary 1.9]{HausenWrobel17}. The structure of the algebra of regular functions on a trinomial hypersurface was studied by Gaifullin and Zaitseva in the papers \cite{GaifullinZaitseva19,Zaitseva19}.

An interesting example of trinomial varieties comes from a construction named suspension. By definition, the suspension over an affine variety $X$ corresponding to a function $f\in\Kp[X]$ is the affine subvariety $\mathrm{Susp}(X,f)$ of $X\times\Ap^2$ defined by the equation $f-uv=0$, where $u,v$ are the coordinate functions on $\Ap^2$. It was proved by Arzhantsev, Zaidenberg and Kuyumzhiyan in~\cite[Theorem 0.2]{ArzhantsevKuyumzhiyanZaidenberg12} that a suspension over a flexible variety is again flexible. In particular, the trinomial hypersurface defined by the equation
$f(T)=T_{21}T_{22}$ is flexible, where $f(T)$ is a polynomial in arbitrary $T_{ij}$ except $T_{21}$ and $T_{22}$. In \cite{Gaifullin19}, Gaifullin proved a sufficient condition for a trinomial hypersurface to be flexible using the correspondence between LNDs and $\Ga$-actions. Precisely, he defined five classes $H_1$--$H_5$ of trinomial hypersurfaces and checked that they are flexible. For instance, a hypersurface of type $H_1$ is defined by the equation
$$T_0^{l_0}+T_1^{l_1}-T_2^1=0,$$
where $T_i^1=T_{i1}T_{i2}\ldots T_{in_i}$. Note that, for $T_2^1=T_{21}T_{22}$ and $f(T)=T_0^{l_0}+T_1^{l_1}$, we obtain the hypersurface $f(T)=T_{21}T_{22}$. For the definition of other types $H_2$--$H_5$ of hypersurfaces, see page~\pageref{def_H} in Section~\ref{sect:main_definitions} below.

The main result of this paper generalizes Gaifullin's results to the case of an arbitrary trinomial variety. Namely, we give a sufficient condition for a trinomial variety to be flexible in terms of degrees of monomials $l_i$ involved in the defining equations of a variety. We also use the correspondence between LNDs and $\Ga$-actions on varieties, but the calculations in the case of an arbitrary trinomial variety become much more technical.

More precisely, we consider five classes of trinomial varieties $V_1$--$V_5$, each of which is obtained from the corresponding trinomial hypersurfaces $H_1$--$H_5$ by adding more equations of the same type. For example, a trinomial variety of type $V_1$ is defined by the system of equations
\begin{equation*}
\begin{cases}
\lambda_2T_0^{l_0} + T_1^{l_1}-T_2^1=0,
\\
\lambda_3T_0^{l_0} + T_1^{l_1}-T_3^1=0,
\\\ldots,\\
\lambda_kT_0^{l_0} + T_1^{l_1}-T_k^1=0.
\end{cases}
\end{equation*}
for some $k\geq3$. Types $V_2$--$V_5$ are defined on page~\pageref{def_V} in the next section. Our main result can be formulated as follows (see Theorem~\ref{thm:flex_V} below).

\medskip\textbf{Theorem.} \emph{Trinomial varieties of types $V_1$\textup, $V_3$ and $V_4$ are flexible}.


\medskip At the contrary, trinomial varieties of types $V_2$ and $V_5$ are not flexible, if they are not hypersurfaces (for type $V_2$, with some additional conditions). Indeed, recall that a variety $X$ is called \emph{rigid} if it does not admit non-trivial $\Ga$-actions. Clearly, if a variety is rigid then it can not be flexible; in some sense, flexibility and rigidity are opposite properties of varieties. In~\cite{EvdokimovaGaifullinShafarevich23}, a criterion of trinomial variety to be rigid in terms of degrees of monomials $l_i$ was given. It follows immediately from this criterion that trinomial varieties of types $V_2$ and $V_5$ are rigid, see Section~\ref{sect:rigid} for the details.

The paper is organized as follows. In Section~\ref{sect:main_definitions}, we briefly recall basic definitions and facts on trinomial varieties and formulate the main result. Section~\ref{Proofs} contains the proofs of certain auxiliary lemmas considering locally nilpotent derivations used in the following. In Section~\ref{sect:flexible} we prove our main result about sufficient conditions of trinomial varieties to be flexible. Section~\ref{sect:rigid} contains some additional results about rigidity of certain type of trinomial varieties.

We thank Ivan Arzhantsev and Sergei Gaifullin for useful discussions.

\sect{The main result: statements} \label{sect:main_definitions}

In this section we give precise definitions and recall basic facts about trinomial varieties and locally nilpotent derivations. After that, we recall a criterium of a trinomial variety to be rigid proved in \cite{EvdokimovaGaifullinShafarevich23}. Finally, we define trinomial varieties of types $V_1$--$V_5$ generalizing Gaifullin's types of hypersurfaces and formulate our main result about flexibility, Theorem~\ref{thm:flex_V}.

We fix integers $r,n > 0, m \geq 0$ and $q \in \{0,1\}$, and a partition 
$$n = n_q + \ldots + n_r,~n_i > 0.$$ 
Let $T_{ij}$ and $S_k$, $q\leq i \leq r,~1 \leq j \leq n$, $1 \leq k \leq m$, be independent variables. We write $\Kp[T_{ij},S_k]$ for the corresponding polynomial ring. For each $i = q,\ldots ,r$, fix a tuple $l_i = (l_{i1},\ldots , l_{in_i})$ of positive integers and define a monomial $$T^{l_i}_i=T^{l_{i1}}_{i1}\ldots T^{l_{in_i}}_{in_i}\in \Kp[T_{ij},~S_k] .$$ Now we introduce the ring $R(A)$ for certain input data~$A$.

Type 1: $q = 1$, $A =(a_1,\ldots , a_n)$, $a_j \in \Kp $ and if $i \neq j$, then $a_i \neq a_j$. Set $I = \{1,\ldots , r-1\}$ and for every $i \in I$ define the polynomial 
$$g_i = T^{l_i}_i - T^{L_{i+1}}_{i+1} - (a_{i+1} - a_i) \in \Kp[T_{ij},S_k].$$
Type 2: $q = 0$, 
$$A=\begin{pmatrix}
  a_{10}& a_{11} & a_{12}  & \ldots & a_{1r}\\
  a_{20}& a_{21} & a_{22}  & \ldots & a_{2r}
\end{pmatrix}$$
is a $2\times (r+1)$ matrix with pairwise linearly independent columns. Set $I = \{0,\ldots , r -2\}$ and define for every $i \in I$ the polynomial
$$g _i = \det \begin{pmatrix}
  T^{l_i}_i& T^{l_{i+1}}_{i+1} & T^{l_{i+2}}_{i+2}\\  
  a_{1i}& a_{1{i+1}} & a_{1{i+2}}\\
  a_{2i}& a_{2{i+1}} & a_{2{i+2}}\\
\end{pmatrix}\in \Kp[T_{ij},S_k].$$
For both types we define $R(A) = \Kp[T_{ij},S_k]/(g_i,i\in I).$ The following definition was given in \cite[Construction 1.1]{HausenWrobel17} (see also \cite{HausenHerppich13,HausenHerppichSuss13,HausenSuss11}).
\defi{\label{def:tr_var} Given a data $A$, the affine variety $$X = \Spec(R(A))$$ is called \emph{trinomial}.
}
It was proved in \cite[Theorem 1.2]{HausenWrobel17} that every trinomial variety is irreducible and normal. Moreover, the following fact was proved.
\mtheo{Suppose $r\geq2$ and $n_il_{ij}>1$ for all $i$\textup, $j$. Then
\textup{(a)} in case of \textup{Type 1}\textup, $R(A)$ is factorial if and only if one has $\mathrm{gcd}(l_{i1},\ldots,l_{in_i}=1$ for all $i=1,\ldots,r$\textup;
\textup{(b)} in case of \textup{Type 2}\textup, $R(A)$ is factorial if and only if the numbers $d_i=\gcd(l_{i1},\ldots,l_{in_i})$ are pairwise coprime.
}

\nota{i) One can easily check that every trinomial variety up to a scalar change of coordinates can be written in the form (\ref{formula:trinomial_variety}), i.e., in the form
\begin{equation*}
X\colon
\begin{cases}
\lambda_2T_0^{l_0} + T_1^{l_1}-T_2^{l_2}=0,\\
\lambda_3T_0^{l_0} + T_1^{l_1}-T_3^{l_3}=0,\\
\ldots,\\
\lambda_kT_0^{l_0} + T_1^{l_1}-T_k^{l_k}=0,
\end{cases}
\end{equation*}
where we keep all the notation from the introduction. So, in the sequel we use this equivalent definition of a trinomial variety.

ii) Of course, if $R(A)$ contains variables $S_k$, then the corresponding trinomial variety is a cylinder over the trinomial variety defined by (\ref{formula:trinomial_variety}) in the affine space with coordinates $T_{ij}$. It is not clear from the notation of (\ref{formula:trinomial_variety}), how many variables $S_k$ are in the ring $R(A)$. We will keep the following convention: everywhere below, except Example~\ref{exam:trinomial_varieties} (ii), we assume that there are no such variables.
}

\exam{i) \label{exam:trinomial_varieties}As it was mentioned in the introduction, $\SL_2(\Kp)$ can be considered as a trinomial hypersurface.

ii) Another interesting class of examples is provided by \emph{Danielewski surfaces}. By definition, such a surface $W_n$ is given by the equation
\begin{equation*}
1+T_{11}^nT_{12}-T_{21}^2=0
\end{equation*}
for certain positive integer $n$ in the affine space with coordinates $T_{11}$, $T_{12}$ and $T_{21}$. As it was shown by W. Danielewski in 1989, these surfaces establish a counterexample to the generalized Zariski cancellation problem, see \cite{Danielewski89} for the details. Namely, he proved that $W_1$ is flexible, while $W_2$ is not, so these surfaces are not isomorphic. But if we consider these subvarieties in the affine space with the additional coordinate $S_1$, then they become isomorphic (in other words, $W_1\times\Ap^1\cong W_2\times\Ap^1$).

iii) The system of equations
\begin{equation}
\begin{cases}
\lambda_2T_{01}^2T_{02}^4 + T_{11}^2T_{12}^6 - T_{21}T_{22}^2=0,\\
\lambda_3T_{01}^2T_{02}^4 + T_{11}^2T_{12}^6-T_{31}T_{32}^5=0\\
\end{cases}\label{formula:TV_V_4}
\end{equation}
defines a 6-dimensional trinomial variety, which is not a trinomial hypersurface in the 8-di\-men\-si\-onal affine space.}

\label{def_H}In \cite{Gaifullin19}, Gaifullin proved a sufficient condition to a trinomial hypersurface to be flexible. Namely, he introduced the following five types of hypersurfaces.
\begin{center}
\begin{tabular}{|l|l|}
\hline
$H_1$&$\Vp(T_0^{l_0} + T_1^{l_1}+T_2^1)\vphantom{\int\limits_a^b}$\\
\hline
$H_2$&$\Vp(T_0^2+T_1^2+T_2^{l_2})\vphantom{\int\limits_a^b}$\\
\hline
$H_3$&$\Vp(T_{01}\wh  T_0^{l_0} + T_1^{l_1}+T_{21}\wh T_2^{l_2})\vphantom{\int\limits_a^b}$\\
\hline
$H_4$&$\Vp(T_{01}^2\wh T_0^{2m_0} + T_{11}^2\wh T_1^{2m_1}+T_{21}\wh T_2^{l_2})\vphantom{\int\limits_a^b}$\\
\hline
$H_5$&$\Vp(T_{01}^2\wh T_0^{2m_0}+T_{11}^2\wh T_1^{2m_1}+T_{21}^2\wh T_2^{2m_2})\vphantom{\int\limits_a^b}$\\
\hline
\end{tabular}
\end{center}
Here we denote 
\begin{equation*}
\begin{split}
T_i^1&=T_{i1}T_{i2}\ldots T_{in_i},~
T_i^2=T_{i1}^2T_{i2}^2\ldots T_{in_i}^2,\\
\wh T_i^{2m_i}&=T_{i2}^{2m_{i2}}T_{i3}^{2m_{i3}}\ldots T_{in_i}^{2m_{in_i}},~
\wh T_i^{l_i}=T_{i2}^{l_{i2}}T_{i3}^{l_{i3}}\ldots T_{in_i}^{l_{in_i}}.
\end{split}
\end{equation*}
As usual, $\Vp(J)$ denotes the set of common zeroes of polynomials from a subset $J\subseteq\Kp[T_{ij}]$. As it was shown in \cite[Theorem 4]{Gaifullin19}, trinomial hypersurfaces of types $H_1$--$H_5$ are flexible.

\exam{The group $\SL_2(\Kp)$ and the Danieliewski surface $W_1$ defined in Section~\ref{sect:main_definitions} are flexible, because they belong to type $H_1$.}

Our goal is to generalize these results. To do this, we introduce the following types of trinomial varieties for an arbitrary $k\geq3$.
\begin{center}
\begin{tabular}{|l|l|}
\hline
$V_1$&$\Vp(\lambda_iT_0^{l_0} + T_1^{l_1}-T_i^1,2\leq i\leq k)\vphantom{\int\limits_a^b}$\\
\hline
$V_2$&$\Vp(\lambda_iT_0^2+T_1^2+T_i^{l_i},2\leq i\leq k)\vphantom{\int\limits_a^b}$\\
\hline
$V_3$&$\Vp(\lambda_iT_{01}\wh T_0^{l_0} + T_1^{l_1}-T_{i1}\wh T_i^{l_i},2\leq i\leq k)\vphantom{\int\limits_a^b}$\\
\hline
$V_4$&$\Vp(\lambda_iT_{01}^2\wh T_0^{2m_0} - T_{11}^2\wh T_1^{2m_1}-T_{i1}\wh T_2^{l_i},2\leq i\leq k)\vphantom{\int\limits_a^b}$\\
\hline
$V_5$&$\Vp(\lambda_iT_{01}^2\wh T_0^{2m_0}+T_{11}^2\wh T_1^{2m_1}+T_{i1}^2\wh T_i^{2m_i},2\leq i\leq k)\vphantom{\int\limits_a^b}$\\
\hline
\end{tabular}
\end{center}
\label{def_V}Our main result can be formulated as follows.
\mtheo{\label{thm:flex_V}Trinomial varieties of types $V_1$\textup, $V_3$ and $V_4$ are flexible.}

At the contrary, a trinomial variety of type $V_5$ is not flexible, while a trinomial variety of type $V_2$ is not flexible under some restriction (we do not know if it is flexible without this restriction, see Section~\ref{sect:rigid} for details).

\sect{Auxiliary lemmas} \label{Proofs}

In this section, we prove three lemmas, which will be used in the next section in the proof of the main result. First, we need to define certain LNDs of the special form.

\cons{\label{rmk:lnds_X}Let
\begin{equation*}
X\colon
\begin{cases}
\lambda_2T_0^{l_0} + T_1^{l_1}-T_2^{l_2}=0,\\
\lambda_3T_0^{l_0} + T_1^{l_1}-T_3^{l_3}=0,\\
\ldots,\\
\lambda_kT_0^{l_0} + T_1^{l_1}-T_k^{l_k}=0.
\end{cases}
\end{equation*}
be a trinomial variety. Suppose that for each $i$ from $2$ to $k$ there exists a number $j_i\in\{1,\ldots,n_i\}$ such that $l_{ij_i}=1$; denote $J=\{j_2,\ldots,j_k\}$. Then for every $p \in \{0,1\}$ and $j_p \in \{1,2,\ldots,n_p\}$ we can define the LNDs $\gamma_{pj_p}^J$ by
\begin{equation*}
\begin{split}
\gamma_{pj_p}^J(T_{pj_p}) &= \prod_{i=2}^k\frac{\dd T_i^{l_i}}{\dd T_{ij_i}},~\gamma_{pj_p}^J(T_{mj_m})
= \frac{\dd T_p^{l_p}}{\dd T_{pj_p}}\prod_{i\in\{2,\ldots k\}\setminus\{m\}}\frac{\dd T_i^{l_i}}{\dd T_{ij_i}} \\
\end{split}
\end{equation*}
for $m=2,\ldots,k$, while $\gamma_{pj_p}^J(T_{ij}) = 0$ for all other pairs of indices $i,j.$ For any $s\in\Kp$, denote $\tau_{pj_p}^J(s) = \exp(s\gamma_{pj_p}^J).$
}

\lemmp{\label{lemm:move_points}Let 
\begin{equation*}
X\colon
\begin{cases}
\lambda_2T_0^{l_0} + T_1^{l_1}-T_2^{l_2}=0,\\
\lambda_3T_0^{l_0} + T_1^{l_1}-T_3^{l_3}=0,\\
\ldots,\\
\lambda_kT_0^{l_0} + T_1^{l_1}-T_k^{l_k}=0,
\end{cases}
\end{equation*}
be a trinomial variety satisfying the conditions of Remark~\textup{\ref{rmk:lnds_X}}. Pick points $P,Q \in X^\reg$. If $\gamma_{ij_i}^J(T_{pj_p})(P) \neq 0$ for some $p \in \{0,1\}$ and $j_p\in\{1,\ldots,n_p\}$ then there exists $s_{ij_i} \in \Kp$ such that $$T_{ij_i}(\tau_{ij_i}^J(s_{ij_i})(P))=T_{ij_i}(Q).$$}
{
We will prove only the case $i =p$, because other cases can be proved similarly.
Since
\begin{equation*}
\tau_{pj_p}^J(s)(T_{pj_p})= T_{pj_p} + s\prod_{i=2}^k\frac{\dd T_i^{l_i}}{\dd T_{ij_i}}\text{ and }\gamma_{pj_p}^J(T_{pj_p})(P)= \prod_{i=2}^k\frac{\dd T_i^{l_i}}{\dd T_{ij_i}} \neq 0,
\end{equation*}
we can put $$s_{pj_p} = \frac{T_{pj_p}(Q)-T_{pj_p}(P)}{\gamma_{pj_p}^J(T_{pj_p})(P)},~R =  \tau_{pj_p}^J(s_{pj_p})(P).$$ Note that, for every $j\in\{1,\ldots,n_p\}\setminus\{j_p\},$ one has $T_{pj_p}(R) = T_{pj_p}(P)$, while $T_{pj_p}(R) = T_{pj_p}(Q)$, as required.
}

\lemmp{\label{lemm:transit_X_C}
Suppose 
\begin{equation*}
X\colon
\begin{cases}
\lambda_2T_0^{l_0} + T_1^{l_1}-T_{21}\wh T_2^{l_2}=0,\\
\lambda_3T_0^{l_0} + T_1^{l_1}-T_{31}\wh T_3^{l_3}=0,\\
\ldots,\\
\lambda_kT_0^{l_0} + T_1^{l_1}-T_{k1}\wh T_k^{l_k}=0.\\
\end{cases}
\end{equation*}
Pick a set of nonzero scalars 
$$C = \{c_{ij}\in\Kp^{\times}\mid2\leq i \leq k,~2\leq j \leq n_i\}.$$ Put $$X_C = \Vp (T_{ij}-c_{ij},2\leq i \leq k,~2\leq j \leq n_i) \cap X.$$ Then $\SAut(X)$ acts on $X_C$ transitively.
}
{
Set $J=\{1,\ldots,1\}$ and recall the notion $\gamma_{pj_p}^1=\gamma_{pj_p}^J$,
Let $P,~Q \in X_C$. Note that

$$\prod_{i=2}^k\wh T_i^{l_i}(P)=\prod_{i=2}^k\wh T_i^{l_i}(Q) \neq 0.$$

By Lemma~\ref{lemm:move_points}, for each $p\in \{0,1\}$ and $j_p=1,\ldots,n_p$, there exists $s_{pj_p}$ such that 
$$T_{pj_p}(R) = T_{pj_p}(Q)\text{ for } R=\tau_{pj_p}^1(s_{pj_p})(P),$$
while if $j\in\{1,\ldots,n_p\}\setminus\{j_p\}$ then $T_{pj_p}(R) = T_{pj_p}(P)$. Starting from~$P$, for every $p \in \{0,1\}$ and $j_p \in \{1,2,\ldots,n_p\}$  we apply the automorphisms $\tau_{pj_p}^1(s_{pj_p})$ step by step to reach a point $S$ such that $T_{pj_p}(S) = T_{pj_p}(Q)$ for every $p$ and $j_p$. Hence,
$$
T_{i1}(S) = \frac{\lambda_iT_0^{l_0}(S)+T_1^{l_1}(S)}{\prod_{i=1}^k\wh T_i^{l_i}(S)} = \frac{\lambda_iT_0^{l_0}(Q)+T_1^{l_1}(Q)}{\prod_{i=1}^k\wh T_i^{l_i}(Q)} = T_{i1}(Q).
$$
Therefore, $S = Q,$ and $P,~Q$ are in the same $\SAut(X)$-orbit, as required.
}

\sect{The main result: proofs} \label{sect:flexible}

This section is devoted to the proof of Theorem~\ref{thm:flex_V}. We will consider the types subsequently. In each type we will check, that the group $\SAut(X)$ acts transitively on regular locus $X^\reg$ of the trinomial variety $X$.

\textit{\textbf{Type}} $V_1$. Suppose $P \in X^\reg$ is such that there exist $p \in \{0,1\}$ and $j_p \in \{1,2,\ldots,n_p\}$ and $\wh J = \{j_3,\ldots,j_k\}$ for which we have 
$$
\prod_{i=3}^k\frac{\dd T_i^{1}}{\dd T_{ij_i}}(P) \neq 0, \frac{\dd T_p^{l_p}}{\dd T_{pj_p}}(P) \neq 0.
$$
Suppose also that a point $Q \in X^\reg$ satisfies $T_{ij}(Q) \neq 0$ for all $i =3,\ldots ,k$ and $j=2,\ldots,n_i.$ By Lemma~\ref{lemm:move_points} for each $m = 2,\ldots,n_2$ there exists $s_{2m}\in\Kp$ such that $\widetilde P = t_{pj_p}^J(s_{2m})(P)$ and $T_{2m}(\widetilde P)=T_{2m}(Q)$. After subsequent applications of $\tau_{pj_p}^J(s_{2m})$ we obtain a point $S$ from the $\SAut(X)$-orbit of $P$ with $T_{2m}(S) = T_{2m}(Q)\neq 0$ for each $m$ from 2 to $n_2$. Since 
$$
\frac{\dd T_2^{1}}{\dd T_{21}}(S) \neq 0,
$$
we can interchange the first and the second equations, put $j_3=1$ and  $\wh J=\{j_3,\ldots,j_k\}$ and build a similar sequence, etc. Then we obtain a point $\widetilde S$ such that $T_{ij}(\widetilde S)=T_{ij}(Q)$ for every $i=2,\ldots,k$ and $j=2,\ldots,n_i.$ By Lemma~\ref{lemm:transit_X_C}, the points $\widetilde S,~Q$ belong to the same $\SAut(X)$-orbit.

Now, pick a point $R\in X^\reg$.  If $T_{pj_p}(R)=0$ for each $p\in\{0,1\}$ and $j_p\in\{1,\ldots,n_p\}$ then there exists $J =\{j_2,\ldots,j_k\}$ such that
$$
\prod_{i=2}^k\frac{\dd T_i^{1}}{\dd T_{ij_i}}(R) \neq 0,
$$
because the point $R$ is regular. Applying, if necessary, the LND $\gamma_{pj_p}^J$ for each $p\in\{0,1\}$ and $j_p\in\{1,\ldots,n_p\}$ and Lemma~\ref{lemm:move_points}, we can assume without loss of generality that $T_{pj_p}\neq 0$ for every $p$ and $j_p$. Then we are in the case when $R,~Q$ are in the same $\SAut(X)$-orbit. 

Finally, suppose that there exist $m$, $t$ such that $T_{m}^{l_m}(R)=T_{t}^{l_t}(R)=0$. Without loss of generality we can assume that $m=2$, $t=3$. Clearly, $T_0^{l_0}(P)=T_1^{l_1}(P)=0$. Since $R\in X^\reg$, there exist $j_0\in\{1,\ldots,n_0\}$ with $l_{0j_0}=1,~j_1\in\{1,\ldots,n_1\}$ with $l_{1j_1}=1$ and $~\wh J=\{j_4,\ldots,j_k\} $ such that
$$
\frac{\dd T_0^{l_0}}{\dd T_{0j_0}}(R)\neq 0,~\frac{\dd T_1^{l_1}}{\dd T_{1j_1}}(R)\neq 0,~\frac{\dd T_i^{1}}{\dd T_{ij_i}}(R)\neq 0.
$$
For $j_3\in\{1,\ldots,n_3\}$, define the LND $\delta_{j_3}$ of $\Kp[X]$ by putting
\begin{equation*}
\begin{split}
\delta_{j_3}(T_{0l_0})&= \frac{\dd T_1^{l_1}}{\dd T_{1j_1}}\prod_{i=3}^k\frac{\dd T_i^{1}}{\dd T_{ij_i}},~
\delta_{j_3}(T_{1l_1})= -\lambda_2\frac{\dd T_0^{l_0}}{\dd T_{0j_0}}\prod_{i=3}^k\frac{\dd T_i^{1}}{\dd T_{ij_i}},\\
\delta_{j_3}(T_{mj_m})&= (\lambda_m-\lambda_2)\frac{\dd T_0^{l_0}}{\dd T_{0j_0}}\frac{\dd T_1^{l_1}}{\dd T_{1j_1}}\prod_{i\neq 2,m}\frac{\dd T_i^{1}}{\dd T_{ij_i}},
\end{split}
\end{equation*}
and $\delta_{j_3}(T_{ij})=0$ for all other pairs of indices $i,j.$ For any $s\in\Kp$ put $\psi_{j_2}(s)=\exp(s\delta_{j_3})$. Since
\begin{equation*}
\begin{split}
&\psi_{j_3}(T_{3j_3}) = T_{3j_3} + s(\lambda_3-\lambda_2)\frac{\dd T_0^{l_0}}{\dd T_{0j_0}}\frac{\dd T_1^{l_1}}{\dd T_{1j_1}}\prod_{i\neq 2,3}\frac{\dd T_i^{1}}{\dd T_{ij_i}},\\
&(\lambda_3-\lambda_2)\frac{\dd T_0^{l_0}}{\dd T_{0j_0}}\frac{\dd T_1^{l_1}}{\dd T_{1j_1}}\prod_{i\neq 2,3}\frac{\dd T_i^{1}}{\dd T_{ij_i}}(R)\neq0,
\end{split}
\end{equation*}
we can find $s_{j_3}\in\Kp$ such that $T_{3j_3}(\psi_{j_3}(R))\neq 0$. By subsequent applying of $\psi_1,\psi_2,\ldots,\psi_{n_3}$ to the point $R$ we reduce the problem to the previous cases.

\textit{\textbf{Type}} $V_3$. Suppose $P\in X^\reg$ is such that $\wh T_{0}^{l_{0}}(P) \neq 0$ and for $i=3,\ldots,k$ one has $\wh T_{i}^{l_{i}}(P)\neq 0$.
Suppose also that a point $Q \in X^\reg$ satisfies $T_{ij}(Q) \neq 0$ for all $i =3,\ldots ,k$ and $j=2,\ldots,n_i$. By Lemma~\ref{lemm:move_points} for each $m = 2,\ldots,n_2$ there exists $s_{2m}\in\Kp$ such that $\widetilde P = t_{pj_p}^J(s_{2m})(P)$ and $T_{2m}(\widetilde P)=T_{2m}(Q)$. After subsequent applications of $\tau_{pj_p}^J(s_{2m})$ we obtain a point $S$ from the $\SAut(X)$-orbit of $P$ with $T_{2m}(S) = T_{2m}(Q)\neq 0$ for each $m$ from 2 to $n_2$. Since 
$$
\frac{\dd T_2^{l_2}}{\dd T_{21}}(S) \neq 0,
$$
we can interchange the first and the second equations, put $j_3=1$ and  $\wh J=\{j_3,\ldots,j_k\}$ and build a similar sequence, etc. Then we obtain a point $\widetilde S$ such that $T_{ij}(\widetilde S)=T_{ij}(Q)$ for every $i=2,\ldots,k$ and $j=2,\ldots,n_i.$ By Lemma~\ref{lemm:transit_X_C}, the points $\widetilde S,~Q$ belong to the same $\SAut(X)$-orbit.

Now, pick a point $R\in X^\reg.$ Assume that $\wh T_{0}^{l_{0}}(R) = 0$ and for each $i=2,\ldots,k$ one has \break$\wh T_{i}^{l_{i}}(R)\neq 0$. We can use LND $\gamma_{0j_0}^1$ and Lemma~\ref{lemm:move_points} to make $\wh T_{0}^{l_{0}}(R)$ be non zero, so we are in the previous case. 
If $\wh T_{0}^{l_{0}}(R) = 0$ and there exists $i=2,\ldots,k$ such that $\wh T_{i}^{l_{i}}(R) = 0$ then either there exist $i$ and $j_i$ such that $l_{ij_i} = 1$ or there exists $j_1$ such that $l_{1j_1} = 1$. In the case $l_{ij_i} = 1$ we can interchange $T_{ij_i}$ and $T_{i0}$ and use a similar line of reasoning. For the case $l_{1j_1} = 1$, we can use $\gamma_{1j_1}^1$ and Lemma~\ref{lemm:move_points} to make $\wh T_{i}^{l_{i}}(R)$ be non zero, so we are again in the previous case.

Finally, suppose that there exist $m$, $t$ such that $T_{m}^{l_m}(R)=T_{t}^{l_t}(R)=0$. Without loss of generality we can assume that $m=2$, $t=3$. Clearly, $T_{00}T_0^{l_0}(P)=T_1^{l_1}(P)=0$. Since $R\in X^\reg$, there exist $j_0\in\{1,\ldots,n_0\}$ with $l_{0j_0}=1,~j_1\in\{1,\ldots,n_1\}$ with $l_{1j_1}=1$ and $~\wh J=\{j_4,\ldots,j_k\} $ such that
$$
\frac{\dd T_0^{l_0}}{\dd T_{0j_0}}(R)\neq 0,~\frac{\dd T_1^{l_1}}{\dd T_{1j_1}}(R)\neq 0,~\frac{\dd T_i^{l_i}}{\dd T_{ij_i}}(R)\neq 0.
$$
For $j_3\in\{1,\ldots,n_3\}$, define the LND $\delta_{j_3}$ of $\Kp[X]$ by
$$
\delta_{j_3}(T_{0l_0}) = \frac{\dd T_1^{l_1}}{\dd T_{1j_1}}\prod_{i=3}^k\frac{\dd T_i^{l_i}}{\dd T_{ij_i}},~ 
\delta_{j_3}(T_{1l_1}) = -\lambda_2\frac{\dd T_0^{l_0}}{\dd T_{0j_0}}\prod_{i=3}^k\frac{\dd T_i^{l_i}}{\dd T_{ij_i}},$$
$$
\delta_{j_3}(T_{mj_m}) = (\lambda_m-\lambda_2)\frac{\dd T_0^{l_0}}{\dd T_{0j_0}}\frac{\dd T_1^{l_1}}{\dd T_{1j_1}}\prod_{i\neq 2,m}\frac{\dd T_i^{l_i}}{\dd T_{ij_i}},$$
and $\delta_{j_3}(T_{ij})=0$ for all other pairs of indices $i,j.$ For any $s\in\Kp$ put $\psi_{j_2}(s)=\exp(s\delta_{j_3})$. Since
\begin{equation*}
\begin{split}
&\psi_{j_3}(T_{3j_3}) = T_{3j_3} + s(\lambda_3-\lambda_2)\frac{\dd T_0^{l_0}}{\dd T_{0j_0}}\frac{\dd T_1^{l_1}}{\dd T_{1j_1}}\prod_{i\neq 2,3}\frac{\dd T_i^{l_i}}{\dd T_{ij_i}},\\
&(\lambda_3-\lambda_2)\frac{\dd T_0^{l_0}}{\dd T_{0j_0}}\frac{\dd T_1^{l_1}}{\dd T_{1j_1}}\prod_{i\neq 2,3}\frac{\dd T_i^{l_i}}{\dd T_{ij_i}}(R)\neq0,
\end{split}
\end{equation*}
we can find $s_{j_3}\in\Kp$ such that $T_{3j_3}(\psi_{j_3}(R))\neq 0$. By subsequent applying of $\psi_1,\psi_2,\ldots,\psi_{n_3}$ to the point $R$ we reduce the problem to the previous cases.

\textit{\textbf{Type}} $V_4$. 
Put 
$$\sqrt{T_0^{2m_0}}=T_{01}\wh T_0^{m_0},~\sqrt{T_1^{2m_1}}=T_{11}\wh T_1^{m_1}.$$ 
Denote $\alpha_i=\sqrt{\lambda_iT_0^{2m_0}}+\sqrt{T_1^{2m_1}}$ and $\beta_i=\sqrt{\lambda_iT_0^{2m_0}}-\sqrt{T_1^{2m_1}}$.
For each vector $J =(j_2,\ldots,j_k)$, $j_i\in\{1,\ldots,n_i\},$ we denote by $\delta_{+}^J$ and $\delta_{-}^J$ two LND's of $\Kp[X]$ given by the formulas
\begin{equation*}
\begin{split}
\delta_{\pm}^J(T_{01})&= \frac{\sqrt{T_1^{2m_1}}}{T_{11}}\prod_{i=2}^k\frac{\dd T_i^{l_i}}{\dd T_{ij_i}} ,~
\delta_{\pm}^J(T_{11})= \mp\frac{\sqrt{T_0^{2m_0}}}{T_{01}}\prod_{i=2}^k\frac{\dd T_i^{l_i}}{\dd T_{ij_i}} ,\\
\delta_{\pm}^J(T_{ij_i})&= (1\pm\sqrt{\lambda_i})\frac{\sqrt{T_0^{2m_0}}}{T_{01}}\frac{\sqrt{T_1^{2m_1}}}{T_{11}}(\sqrt{T_0^{2m_0}}\pm \sqrt{\lambda_i}\sqrt{T_1^{2m_1}})\\
&\times\prod_{s\neq i}\frac{\dd T_s^{l_s}}{\dd T_{sj_s}}\text{ for all }i=2,\ldots,k,\\
\end{split}
\end{equation*}
and $\delta_+^J(T_{ij})=\delta_-^J(T_{ij})=0$ for all other $i,j.$ 

Suppose $Q\in X^\reg$ is such that, for $i=3,\ldots,k$, one has $\wh T_{i}^{l_{i}}(Q)\neq 0$ and $\alpha_i(Q)\neq 0$. Note that if there exists $i=2,\ldots,k$ for which $\alpha_i(Q)\neq0$ or $\beta_i(Q)\neq0$, then for all $i=2,\ldots,k$ we have $\alpha_i(Q)\neq0$ or $\beta_i(Q)\neq0$. Pick also a point $P\in X^\reg$ such that, for $i=2,\ldots,k$, one has $\wh T_{i}^{l_{i}}(P)\neq 0.$ Now, we can use the exponents of the LND's defined in Remark~\ref{rmk:lnds_X} to move the point $P$ to a point $\widetilde P\in X^\reg$ such that $\alpha_i(\widetilde P)=\alpha_i(Q)\neq0$ for all $i$ from 2 to $k$. Then for each $j_2=2,\ldots,n_i$, using a composition of $\exp(s_{j_2}\delta_+^J)$, where $J=(j_2,1,\ldots,1)$, we can obtain a point $R\in X^\reg$ such that $T_{2j_2}(P)=T_{2j_2}(Q)$. Those, we are in the case of Lemma~\ref{lemm:transit_X_C}, so $P$ and $Q$ are in the same $\SAut(X)$-orbit, as required.

Next, suppose that for a point $P\in X^\reg$ there exists $i=2,\ldots,k$, such that $\wh T_{i}^{l_{i}}(P)=0.$ We may assume without loss of generality that $i=2$. Then, since $P$ is a regular point of $X$, at least one of the following conditions is satisfied: there exists $j_2=1,\ldots,n_2$ such that $l_{2j_2}=1$ and $\dfrac{\dd T_2^{l_2}}{\dd T_{2j_2}}(P)\neq 0$, or $\alpha_2(P)\neq0$, or $\beta_2(P)\neq0$. In the first case we can interchange $T_{2j_2}$ and $T_{21}$. Using a reasoning same to the previous paragraph, we  obtain that the points $P$ and $Q$ are in the same $\SAut(X)$-orbit. The second and the third cases can be proved completely similarly, so we will prove only the case with $\alpha_i(P)\neq 0$. Here, for each $j_2=2,\ldots,n_i$, using a composition of $\exp(s_{j_2}\delta_+^J)$, where $J=(j_2,1,\ldots,1)$, we obtain a point $R\in X^\reg$ such that $T_{2j_2}(R)\neq0$, and thus we are in the previous case.


\exam{\label{ex:V1-V5}i) The variety
\begin{equation*}
X\colon
\begin{cases}
T_{01}^5T_{02}^4+T_{11}T_{12}^7-T_{21}T_{22}^8T_{23}=0,\\
2T_{01}^5T_{02}^4+T_{11}T_{12}^7-T_{31}T_{32}^2T_{33}^4=0,
\end{cases}
\end{equation*}
is flexible, because it is of type $V_3$.

ii) The variety defined by formula (\ref{formula:TV_V_4}) in Example~\ref{exam:trinomial_varieties} (iii) is flexible, because it belongs to type $V_4$.

iii) The trinomial hypersurface
$$T_{01}^3+T_{11}^5+T_{21}T_{22}T_{23}^2=0,$$
as well as the trinomial variety
\begin{equation}
X\colon
\begin{cases}
\lambda_2T_0^2+T_1^2-T_2^{l_2}=0,\\
\lambda_3T_0^2+T_1^2-T_{31}\wh T_3^{l_3}=0,\\
\ldots,\\
\lambda_kT_0^2+T_1^2-T_{k1}\wh T_k^{l_k}=0,
\end{cases}
\label{ex:unknown}
\end{equation}
where the tuple $l_2$ does not contain $1$,
do not belong to the types under consideration, and, in fact, we do not know if they are flexible or not.}

\sect{Concluding remarks} \label{sect:rigid}

It was proved in the previous section that trinomial varieties of types $V_1$, $V_3$ and $V_4$ are flexible. At the contrary, trinomial varieties of type $V_5$ are not flexible, while trinomial varieties of type $V_2$ are not flexible under some restriction, and we can not prove that they are not flexible without this restriction. To prove this, we recall the notion of rigidity: a variety is called rigid if it does not admit non-trivial $\Ga$-actions. In \cite{EvdokimovaGaifullinShafarevich23}, for a trinomial variety, a criterion to be rigid was proved.

\mtheo{\textup{\cite[Theorem 1]{EvdokimovaGaifullinShafarevich23}}
Let $X$ be a trinomial variety of Type \textup{1}. Then $X$ is not rigid if and only if one of the following holds\textup{:}
\newline\textup{1)} $m > 0$\textup;
\newline\textup{2)} there is $b \in \{1,\ldots,r\}$ such that for each $$i \in \{1,\ldots,r\}\setminus\{b\}$$ there is $j(i)\in\{1,\ldots,n_i\}$ with $l_{ij(i)} = 1.$
}
\mtheo{\label{thm:rig_t2}\textup{\cite[Theorem 3]{EvdokimovaGaifullinShafarevich23}}
Let $X$ be a trinomial variety of Type \textup2. Then $X$ is not rigid if and only if one of the following holds\textup{:}
\newline\textup{1)} $m > 0$\textup{;}
\newline\textup{2)} there are at most two numbers $a,b \in \{0,\ldots,r\}$ such that for each $$i \in \{0,\ldots,r\}\setminus\{a,b\}$$ there is $j(i)\in\{1,\ldots,n_i\}$ such that $l_{ij(i)} = 1$\textup{;}
\newline\textup{3)} there are exactly three numbers $a,b,c\in \{0,\ldots,r\}$ such that for each $i \in \{a,b\}$ there is\break $j(i)\in \{1,\ldots,n_i\}$ with $l_{ij(i)}=2$ and the numbers $l_{ik}$ are even for all $k\in\{1,\ldots,n_i\}$. Moreover\textup{,} for each $i \in \{0,\ldots,r\}\setminus\{a,b,c\}$\textup{,} there is $j(i) \in \{1,\ldots,n_i\}$ with $l_{ij(i)} = 1.$
}
A trinomial variety of type $V_5$ belongs to type $2$. One can check that such a variety does not satisfy the conditions of Theorem~\ref{thm:rig_t2}. The entire set of monomials is
\begin{equation*}
T=\{T_{01}^2\wh T_0^{2m_0},T_{11}^2\wh T_1^{2m_1},\\
T_{21}^2\wh T_2^{2m_2},\ldots,T_{k1}^2\wh T_k^{2m_k}\}.
\end{equation*}
Hence, such variety is rigid, and, consequently, is not flexible.

On the other hand, for a trinomial variety of type $V_2$, the entire set of monomials is
$$T=\{T_0^2,T_1^2,T_2^{l_2}\ldots,T_k^{l_k}\}.$$
Such a variety also belongs to type 2. Assume that there exist at least two indices $i_1,i_2\in\{2,\ldots,k\}$, whose monomials do not contain any variable of degree $1$, i.e., there are no $$j(i_s)\in\{1,\ldots,\ n_{i_s}\},~s=1,2,$$
for which $l_{i_sj(i_s)}=1$. In this case, it is easy to see that our variety does not satisfy the conditions of Theorem~\ref{thm:rig_t2}.

Finally, the trinomial variety, defined in Example~\ref{ex:V1-V5} (iv) by formula (\ref{ex:unknown}) satisfies the conditions of Theorem~\ref{thm:rig_t2}, so it is not rigid. At the moment, we can neither prove that it is flexible nor give a counterexample.

\bigskip\textsc{Mikhail Ignatev: National Research University Higher School of Economics,\break Pokrovsky Boulevard 11, 109028, Moscow, Russia}

\emph{E-mail address}: \texttt{mihail.ignatev@gmail.com}

\bigskip\textsc{Timofey Vilkin: National Research University Higher School of Economics,\break Pokrovsky Boulevard 11, 109028, Moscow, Russia}

\emph{E-mail address}: \texttt{vilkin.timofey@gmail.com}

\end{document}